\PassOptionsToPackage{unicode}{hyperref}
\PassOptionsToPackage{hyphens}{url}
\PassOptionsToPackage{dvipsnames,svgnames,x11names}{xcolor}
\documentclass[
]{article}

\usepackage{amsmath,amssymb}
\usepackage{iftex}
\ifPDFTeX
  \usepackage[T1]{fontenc}
  \usepackage[utf8]{inputenc}
  \usepackage{textcomp} 
\else 
  \usepackage{unicode-math}
  \defaultfontfeatures{Scale=MatchLowercase}
  \defaultfontfeatures[\rmfamily]{Ligatures=TeX,Scale=1}
\fi
\usepackage{lmodern}
\ifPDFTeX\else  
\fi
\IfFileExists{upquote.sty}{\usepackage{upquote}}{}
\IfFileExists{microtype.sty}{
  \usepackage[]{microtype}
  \UseMicrotypeSet[protrusion]{basicmath} 
}{}
\makeatletter
\@ifundefined{KOMAClassName}{
  \IfFileExists{parskip.sty}{%
    \usepackage{parskip}
  }{
    \setlength{\parindent}{0pt}
    \setlength{\parskip}{6pt plus 2pt minus 1pt}}
}{
  \KOMAoptions{parskip=half}}
\makeatother
\usepackage{xcolor}
\makeatletter
\ifx\paragraph\undefined\else
  \let\oldparagraph\paragraph
  \renewcommand{\paragraph}{
    \@ifstar
      \xxxParagraphStar
      \xxxParagraphNoStar
  }
  \newcommand{\xxxParagraphStar}[1]{\oldparagraph*{#1}\mbox{}}
  \newcommand{\xxxParagraphNoStar}[1]{\oldparagraph{#1}\mbox{}}
\fi
\ifx\subparagraph\undefined\else
  \let\oldsubparagraph\subparagraph
  \renewcommand{\subparagraph}{
    \@ifstar
      \xxxSubParagraphStar
      \xxxSubParagraphNoStar
  }
  \newcommand{\xxxSubParagraphStar}[1]{\oldsubparagraph*{#1}\mbox{}}
  \newcommand{\xxxSubParagraphNoStar}[1]{\oldsubparagraph{#1}\mbox{}}
\fi
\makeatother

\usepackage{longtable,booktabs,array}
\usepackage{calc} 
\usepackage{etoolbox}
\makeatletter
\patchcmd\longtable{\par}{\if@noskipsec\mbox{}\fi\par}{}{}
\makeatother
\IfFileExists{footnotehyper.sty}{\usepackage{footnotehyper}}{\usepackage{footnote}}
\makesavenoteenv{longtable}
\usepackage{graphicx}
\makeatletter
\def\maxwidth{\ifdim\Gin@nat@width>\linewidth\linewidth\else\Gin@nat@width\fi}
\def\maxheight{\ifdim\Gin@nat@height>\textheight\textheight\else\Gin@nat@height\fi}
\makeatother
\setkeys{Gin}{width=\maxwidth,height=\maxheight,keepaspectratio}
\makeatletter
\def\fps@figure{htbp}
\makeatother

\usepackage{amsmath}
\usepackage{tikz}
\usepackage{xcolor}
\usepackage{arxiv}
\usepackage{orcidlink}
\usepackage{amsmath}
\usepackage[T1]{fontenc}
\makeatletter
\@ifpackageloaded{caption}{}{\usepackage{caption}}
\AtBeginDocument{%
\ifdefined\contentsname
  \renewcommand*\contentsname{Table of contents}
\else
  \newcommand\contentsname{Table of contents}
\fi
\ifdefined\listfigurename
  \renewcommand*\listfigurename{List of Figures}
\else
  \newcommand\listfigurename{List of Figures}
\fi
\ifdefined\listtablename
  \renewcommand*\listtablename{List of Tables}
\else
  \newcommand\listtablename{List of Tables}
\fi
\ifdefined\figurename
  \renewcommand*\figurename{Figure}
\else
  \newcommand\figurename{Figure}
\fi
\ifdefined\tablename
  \renewcommand*\tablename{Table}
\else
  \newcommand\tablename{Table}
\fi
}
\@ifpackageloaded{float}{}{\usepackage{float}}
\floatstyle{ruled}
\@ifundefined{c@chapter}{\newfloat{codelisting}{h}{lop}}{\newfloat{codelisting}{h}{lop}[chapter]}
\floatname{codelisting}{Listing}

\makeatother
\makeatletter
\@ifpackageloaded{caption}{}{\usepackage{caption}}
\@ifpackageloaded{subcaption}{}{\usepackage{subcaption}}
\makeatother

\ifLuaTeX
  \usepackage{selnolig}  
\fi
\usepackage[numbers]{natbib}
\usepackage{bookmark}

\IfFileExists{xurl.sty}{\usepackage{xurl}}{} 
\hypersetup{
  pdftitle={KKT-Informed Neural Network},
  pdfauthor={Carmine Delle Femine},
  pdfkeywords={Optimization, Parametric Optimization, Convex
Optimization, Karush-Kuhn-Tucker (KKT) Conditions, Neural Networks},
  colorlinks=true,
  linkcolor={blue},
  filecolor={Maroon},
  citecolor={Blue},
  urlcolor={Blue},
  pdfcreator={LaTeX via pandoc}}

\renewcommand{\today}{September 11, 2024}
\newcommand{\runninghead}{A Preprint }
\title{KKT-Informed Neural Network}
\usepackage{etoolbox}
\makeatletter
\providecommand{\subtitle}[1]{
  \apptocmd{\@title}{\par {\large #1 \par}}{}{}
}
\makeatother
\subtitle{A Parallel Solver for Parametric Convex Optimization Problems}
\author{\textbf{Carmine Delle Femine}\\Department of Data Intelligence
for Energy and Industrial Processes\\Vicomtech Foundation\\Donostia-San
Sebastián,\ 20009\\\href{mailto:cdellefemine@vicomtech.org}{cdellefemine@vicomtech.org}}
\date{September 11, 2024}
\begin{document}
\maketitle
\begin{abstract}
A neural network-based approach for solving parametric convex
optimization problems is presented, where the network estimates the
optimal points given a batch of input parameters. The network is trained
by penalizing violations of the Karush-Kuhn-Tucker (KKT) conditions,
ensuring that its predictions adhere to these optimality criteria.
Additionally, since the bounds of the parameter space are known,
training batches can be randomly generated without requiring external
data. This method trades guaranteed optimality for significant
improvements in speed, enabling parallel solving of a class of
optimization problems.
\end{abstract}
{\bfseries \emph Keywords}
\def\sep{\textbullet\ }
Optimization \sep Parametric Optimization \sep Convex
Optimization \sep Karush-Kuhn-Tucker (KKT) Conditions \sep 
Neural Networks

\section{Introduction}\label{introduction}

Solving convex optimization problems is essential across numerous
fields, including optimal control, logistics, and finance. In many
scenarios, such as the development of surrogate models, there is a need
to solve a large set of related optimization problems defined by varying
parameters. Achieving fast solutions, even at the cost of strict
optimality guarantees, is often a priority.

Neural networks, with their inherent ability to process data in parallel
and adapt to diverse problem structures, offer a promising solution.
This work introduces the KKT-Informed Neural Network (KINN), a method
designed to solve parametric convex optimization problems efficiently by
integrating the KKT conditions into the network's learning process. This
approach enables rapid, parallel problem-solving while balancing the
trade-off between speed and guaranteed optimality.

\section{Background}\label{background}

Consider a parametric convex optimization problem in the standard form:

\[
\begin{aligned}
\min_{x \in \mathcal{D} \subseteq\mathbb{R}^n} \quad &f(x, {\theta})\\
\textrm{s.t.} \quad & g_i(x, \theta) \leq 0 \quad i = 1, \dots, m \\
& A(\theta) x - b(\theta) = 0
\end{aligned}
\]

where \(x \in \mathcal{D} \subseteq\mathbb{R}^n\) is the optimization
variable; \(\theta \in \mathcal{D}_\theta \subseteq \mathbb{R}^k\) are
the parameters defining the problem;
\(f: \mathcal{D}_f \subseteq\mathbb{R}^n \times \mathbb{R}^k \to \mathbb{R}\)
is the convex cost function;
\(g_i: \mathcal{D}_{g_i} \subseteq\mathbb{R}^n \times \mathbb{R}^k \to \mathbb{R}\)
are the convex inequality constraints,
\(A: \mathcal{D}_\theta \to \mathbb{R}^{p \times n}\) and
\(b: \mathcal{D}_\theta \to \mathbb{R}^{p}\) defines the affine equality
constraints and
\(\mathcal{D} = \bigcap_{i=1}^{m} \mathcal{D}_{g_i} \cap \mathcal{D}_{f}\)
is the domain of the optimization problem.

Assume differentiable cost and constraints functions and that \(g_i\)
satisfies Slater's condition. Given a set of parameters \(\theta\),
\(x^* \in \mathcal{D}\) is optimal if and only if there are
\(\lambda^*\) and \(\nu^*\) that, with \(x^*\), satisy the
Karush-Kuhn-Tucker conditions (KKT) \citep{boydConvexOptimization2004}:

\begin{align}
    A(\theta) x^* - b(\theta) = 0&\\
    g_i(x^*, \theta) \leq 0& \quad i=1,\dots, m\\
    \lambda_i^* \geq 0& \quad i=1,\dots, m\\
    \lambda_i^* g_i(x^*, \theta) = 0& \quad i=1,\dots, m\\
    \nabla_{x^*} f(x^*, \theta) + \sum\nolimits_{i=1}^m \lambda^*_i\nabla_{x^*} g_i(x^*, \theta) + A(\theta)^T\nu^* = 0 &
\end{align}

\section{Proposed method}\label{proposed-method}

KKT-Informed Neural Network (KINN) builds upon the principles of
Physics-Informed Neural Networks (PINNs)
\citep{raissiPhysicsinformedNeuralNetworks2019a}, inducing compliance
with Karush-Kuhn-Tucker (KKT) through a learning bias, directly coding
their violation into the loss function that will be minimized in the
training phase.

The network is designed as a multi-layer perceptron (MLP) and processes
a batch of \(B\) problem parameters
\(\Theta \in \mathbb{R}^{B \times k}\), \(\Theta_i = \theta^{(i)}\). The
network outputs \(\hat{X}\), \(\hat\Lambda\), \(\hat{N}\). A ReLU
function is applied to the branch predicting \(\hat\Lambda\) to ensure
its feasibility.

\begin{align}
[\hat{X}, \hat{\Lambda}, \hat{N}] &= \textrm{KINN}(\Theta)\\
\hat{X} \in \mathbb{R}^{B\times n}&, \quad \hat X_i = \hat{x}^{(i)}\\
\hat{\Lambda} \in \mathbb{R}^{0^{B\times m}}_+&, \quad \hat\Lambda_i = \hat\lambda^{(i)}\\
\hat{N} \in \mathbb{R}^{B\times p}&, \quad \hat{N}_i = \hat\nu^{(i)}
\end{align}

Vector-valued loss function consists of four terms that correspond to
each KKT conditions:

\begin{equation}
\mathcal{L} = \frac{1}{B}\biggr[\sum_{i=1}^B\mathcal{L}_{S}^{(i)}, \sum_{i=1}^B\mathcal{L}_{I}^{(i)}, \sum_{i=1}^B\mathcal{L}_{E}^{(i)}, \sum_{i=1}^B\mathcal{L}_{C}^{(i)}\biggr] 
\end{equation}

where:

\begin{align}
    \mathcal{L}_{S}^{(i)} =& \|\nabla_{\hat{x}^{(i)}} f(\hat{x}^{(i)}, \theta^{(i)}) + \sum\nolimits_{j=1}^m \hat{\lambda}^{(i)}_j\nabla_{\hat{x}^{(i)}} g_j(\hat{x}^{(i)}, \theta^{(i)}) + A(\theta^{(i)})^T\hat{\nu}^{(i)}\|_2\\ 
    \mathcal{L}_{I}^{(i)}  =& \|(\max(0, g_1(\hat{x}^{(i)}, \theta^{(i)})),\dots,\max(0, g_m(\hat{x}^{(i)}, \theta^{(i)})))\|_2\\
    \mathcal{L}_{E}^{(i)} =& \|A(\theta^{(i)}) \hat{x}^{(i)} - b(\theta^{(i)})\|_2\\
    \mathcal{L}_{C}^{(i)}  =& \|(\hat{\lambda}_1^{(i)} g_1(\hat{x}^{(i)}, \theta^{(i)}),\dots,\hat{\lambda}_m^{(i)} g_m(\hat{x}^{(i)}, \theta^{(i)}))\|_2\\
\end{align}

This vector-valued loss function is minimized through a Jacobian descent
\citep{quintonJacobianDescentMultiObjective2024a}. Let
\(\mathcal{J} \in \mathbb{R}^{P\times4}\) the Jacobian matrix of
\(\mathcal{L}\), with \(P\) the number of parameters of the KINN,
\(\mathcal{A}: \mathbb{R}^{P\times4} \to \mathbb{R}^{P}\) is called
aggregator. The ``direction'' of the update of networks parameters will
be \(\mathcal{A}(\mathcal{J})\). The aggregator chosen is
\(\mathcal{A}_{\mathrm{UPGrad}}\), described in
\citep{quintonJacobianDescentMultiObjective2024a}.

\section{Case study}\label{case-study}

A renewable energy generator in a power grid is used as a test case for
this approach.. The generator's active and reactive power injections
\((P,Q)\) are controllable, but they must adhere to physical
constraints. As such, the desired setpoints \((a_P, a_Q)\) must be
projected onto the feasible set defined by these constraints. This
problem is taken from \citep{henryGymANMReinforcementLearning2021}.

\subsection{Problem description}\label{problem-description}

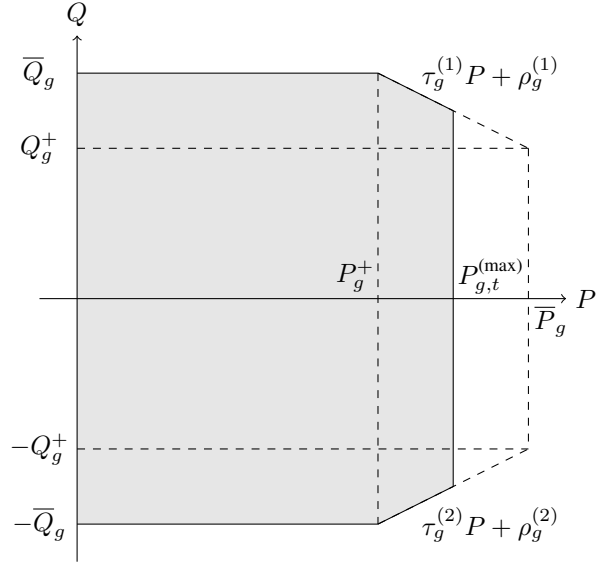
\begin{figure}
\begin{center}
\begin{tikzpicture}


\fill[gray!20] (0,-3) rectangle (4,3);
\fill[gray!20] (4,-2) rectangle (5,2);
\fill[gray!20] (4,2) -- (5,2) -- (5, 2.5) -- (4,3)-- cycle;
\fill[gray!20] (4,-2) -- (5,-2) -- (5, -2.5) -- (4,-3) -- cycle;
\draw (0,3) -- (4,3);
\draw[dashed] (0,2) -- (6,2);
\draw[dashed] (0,-2) -- (6,-2);
\draw (0,-3) -- (4,-3);

\draw[dashed] (4,3) -- (4,-3);
\draw (5,2.5) -- (5,-2.5);
\node at (-0.5,3) {$\overline{Q}_g$};
\node at (-0.5,2) {$Q_g^+$};
\node at (-0.5,-2) {$-Q_g^+$};
\node at (-0.5,-3) {$-\overline{Q}_g$};

\node at (3.7,0.3) {$P_g^+$};
\node at (5.5,0.3) {$P^{\textrm{(max)}}_{g,t}$};
\node at (6.3,-0.3) {$\overline{P}_g$};

\draw (4,3) -- (5,2.5);
\draw (4,-3) -- (5,-2.5);
\draw[dashed] (4,3) -- (6,2);
\draw[dashed] (4,-3) -- (6,-2);

\draw[dashed] (6,2) -- (6,-2);
\node at (5.5,3) {$\tau_g^{(1)} P + \rho_g^{(1)}$};
\node at (5.5,-3) {$\tau_g^{(2)} P + \rho_g^{(2)}$};

\draw[->] (-0.5,0) -- (6.5,0) node[right] {$P$};
\draw[->] (0,-3.5) -- (0,3.5) node[above] {$Q$};
\end{tikzpicture}
\end{center}
\caption{Feasibile set $\mathcal{D}$} \label{fig:set}
\end{figure}

The feasible set \(\mathcal{D}\) (shown in Figure \(\ref{fig:set}\)) is
defined by the physical parameters of the generator
\(\overline{P}_g \in \mathbb{R}_0^+\),
\(P^+_g \in ]0, \overline{P}_g]\),
\(\overline{Q}_g \in \mathbb{R}_0^+\),
\(Q^+_g \in ]0, \overline{Q}_g]\), characterizing the minimum and
maximum possible values and the relationships between active and
reactive power, and the dynamic value \(P^{\textrm{(max)}}_{g,t}\) which
indicates the maximum power that can be generated at that time given the
external conditions (e.g.~wind speed, solar radiation, etc.):
\begin{equation}
\mathcal{D} = \{(P, Q) \in \mathbb{R}² | 0 \leq P \leq P^{\textrm{(max)}}_{g,t}, -\overline{Q}_g \leq Q \leq \overline{Q}_g, Q \leq \tau^{(1)}_g P + \rho_g^{(1)}, Q \geq \tau^{(2)}_g P + \rho_g^{(2)}\}
\end{equation}

where:

\begin{align}
    \tau^{(1)}_g &= \frac{Q_g^+ - \overline{Q}_g}{\overline{P_g} - P_g^+}\\
    \rho^{(1)}_g &= \overline{Q}_g - \tau^{(1)}_gP_g^+\\
    \tau^{(2)}_g &= \frac{\overline{Q}_g -Q_g^+ }{\overline{P_g} - P_g^+}\\
    \rho^{(2)}_g &= -\overline{Q}_g - \tau^{(2)}_gP_g^+  \\
\end{align}

The problem could be stated in standard form as:

\[
\begin{aligned}
\min_{x \in \mathcal{D} \subseteq\mathbb{R}^2} \quad &\frac{1}{2}\|a - x\|^{2}_{2} \\
\textrm{s.t.}\quad & G x - h \leq 0\\
\end{aligned}
\]

with \(a = (a_P, a_Q)\), \(x = (P, Q)\) and:

\begin{equation}
G = \begin{pmatrix}
-1 & 1 & 1 & 0 & 0 & -\tau^{(1)}_g & \tau^{(2)}_g\\
0 & 0 & 0 & -1 & 1 & 1 & 1
\end{pmatrix}^T
\end{equation}

\begin{equation}
h = \begin{pmatrix}
0 & \overline{P}_g & P^{\textrm{(max)}}_{g,t} & \overline{Q}_g & \overline{Q}_g & \rho^{(1)}_g & -\rho^{(2)}_g
\end{pmatrix}^T
\end{equation}

With associated KKT conditions: \begin{align}
    G{x}^* - h \leq 0&\\
    \lambda_i^* \geq 0& \quad i=1,\dots, 7\\
    G^T\lambda^* = 0&\\
    (a-{x}^*) + G^T\lambda^*  = 0 &
\end{align}

\subsection{Experimental results}\label{experimental-results}

The problem described has a two-dimensional optimization variable, seven
scalar parameters and seven constraints:

\begin{equation}
 [\hat{X}, \hat\Lambda] = \mathrm{KINN}(\Theta)
\end{equation}

with: \begin{align}
\Theta \in \mathbb{R}^{B \times 7}, \quad&\Theta_{i} = (a_P^{(i)}, a_Q^{(i)},\overline{P}_g^{(i)},  P_g^{+^{(i)}}, \overline{Q}_g^{(i)},  Q_g^{+^{(i)}}, P^{\textrm{(max)}^{(i)}}_{g,t})\\
\hat{X} \in \mathbb{R}^{B \times 2}, \quad&\hat{X}_{i} = \hat{x}^{(i)} = (\hat{P}^{(i)}, \hat{Q}^{(i)})\\
\hat\Lambda \in \mathbb{R}_+^{0^{B \times 7}}, \quad&\hat\Lambda_{i} = \hat\lambda^{(i)}
\end{align}

The network is composed by three hidden layers of \(512\) neurons each,
with a LeakyReLU (negative slope of \(0.01\)) as activation function and
a skip connection around each hidden layer.

At each training step, a random batch of parameters \(\Theta\) was
sampled:

\begin{align}
a_P^{(i)} &\sim U(0~\mathrm{p.u}., 1~\mathrm{p.u.})\\
a_Q^{(i)} &\sim U(-1~\mathrm{p.u}., 1~\mathrm{p.u.})\\
\overline{P}_g^{(i)} &\sim U(0.2~\mathrm{p.u}., 0.8~\mathrm{p.u.})\\
P_g^{+^{(i)}} &\sim U(0~\mathrm{p.u}., \overline{P}_g^{(i)})\\
\overline{Q}_g^{(i)} &\sim U(0.2~\mathrm{p.u}., 0.8~\mathrm{p.u.})\\
Q_g^{+^{(i)}} &\sim U(0~\mathrm{p.u}., \overline{Q}_g^{(i)})\\
P^{\textrm{(max)}^{(i)}}_{g,t} &\sim U(0~\mathrm{p.u}., \overline{P}_g^{(i)})
\end{align}

Models parameters were update to minimize the following vector-valued
loss function:

\begin{equation}
\mathcal{L} = \frac{1}{B}\biggr[\sum_{i=1}^B\mathcal{L}_{S}^{(i)}, \sum_{i=1}^B\mathcal{L}_{I}^{(i)}, \sum_{i=1}^B\mathcal{L}_{C}^{(i)}\biggr] 
\end{equation}

where:

\begin{align}
    \mathcal{L}_S^{(i)} =& \|(a^{(i)}-\hat{x}^{(i)}) + G^{(i)^{T}}\hat\lambda^{(i)}\|_2\\ 
    \mathcal{L}_{I}^{(i)}  =& \|\max(0, G^{(i)}\hat{x} - h^{(i)})\|_2\\
    \mathcal{L}_{C}^{(i)}  =& \|G^{(i)^{T}}\hat\lambda^{(i)}\|_2\\
\end{align}

\subsubsection{Training}\label{training}

Training was performed with the Adam optimization algorithm with an
intial learning rate of \(10^{-3}\) and an exponential scheduler for the
latter with a \(\gamma\) of \(0.99986\). An early stopping condition
occurred when no progress occurred on any of the constituent terms of
the loss for \(5000\) steps.

\begin{figure}

\centering{

\includegraphics{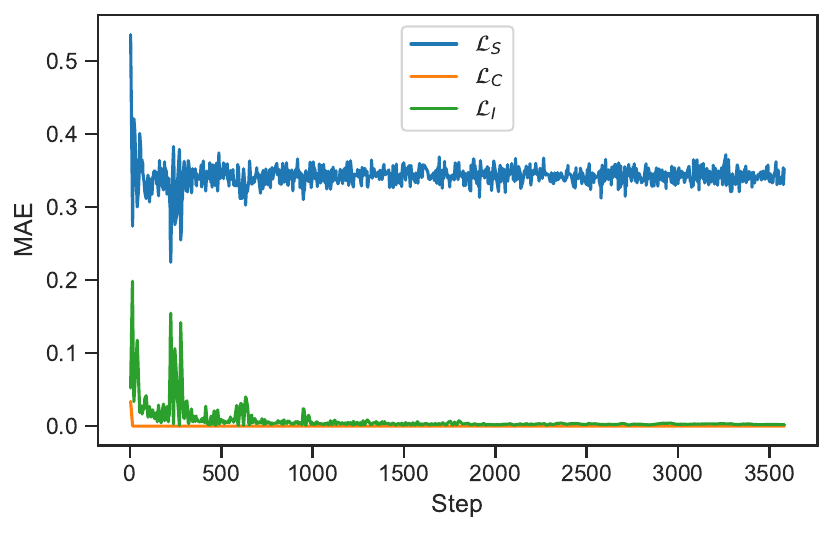}

}

\caption{\label{fig-training}Loss terms during training}

\end{figure}

Finally, the training lasted \(3583\) steps reaching final values shown
in Table~\ref{tbl-training}, while the evolution along the various steps
is in Figure~\ref{fig-training}.

\begin{longtable}[]{@{}cc@{}}
\caption{Final values}\label{tbl-training}\tabularnewline
\toprule\noalign{}
Loss & Value \\
\midrule\noalign{}
\endfirsthead
\toprule\noalign{}
Loss & Value \\
\midrule\noalign{}
\endhead
\bottomrule\noalign{}
\endlastfoot
\(\mathcal{L}_S\) & 0.3519 \\
\(\mathcal{L}_I\) & 0.0020 \\
\(\mathcal{L}_C\) & 0.0000 \\
\end{longtable}

\subsubsection{Evaluation}\label{evaluation}

To evaluate the approach presented here, the ``cvxpylayers'' library
\citep{agrawalDifferentiableConvexOptimization2019a}, which implements a
batched solver via multi-threading, was used as a baseline. The
validation set consists of 1000 samples, generated by taking the
physical parameters of the two generators present in the use case
presented in \citep{henryGymANMReinforcementLearning2021} and, for each
of them, simulating 500 random inputs \((a_P, a_Q)\) and external
condition \(P^{\textrm{(max)}}_{g,t}\). Specifically these parameters
are:

\begin{align}
    [\overline{P}_g^{(i)},  P_g^{+^{(i)}}, \overline{Q}_g^{(i)},  Q_g^{+^{(i)}}]_1 &= [0.3, 0.2, 0.3, 0.15]\\
    [\overline{P}_g^{(i)},  P_g^{+^{(i)}}, \overline{Q}_g^{(i)},  Q_g^{+^{(i)}}]_2 &= [0.5, 0.35, 0.5, 0.2]
\end{align}

\begin{figure}

\begin{minipage}{0.50\linewidth}

\centering{

\includegraphics{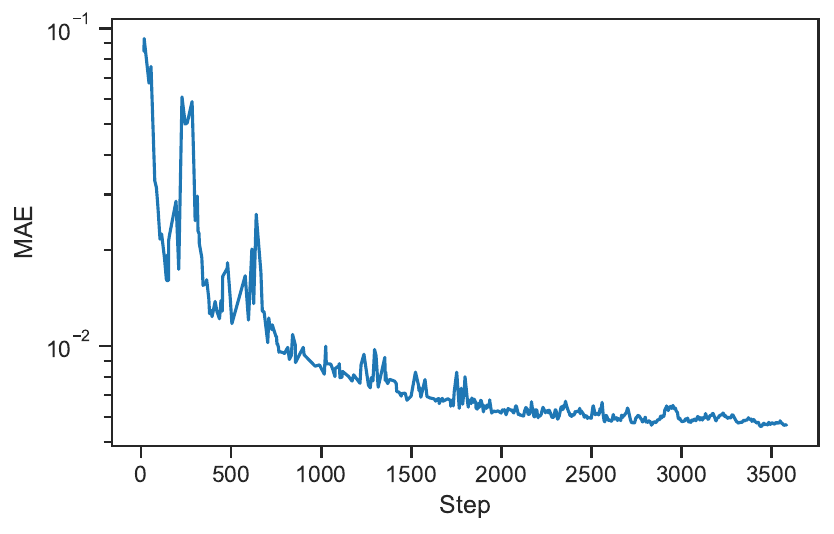}

}

\subcaption{\label{fig-val-1}MAE}

\end{minipage}%
\begin{minipage}{0.50\linewidth}

\centering{

\includegraphics{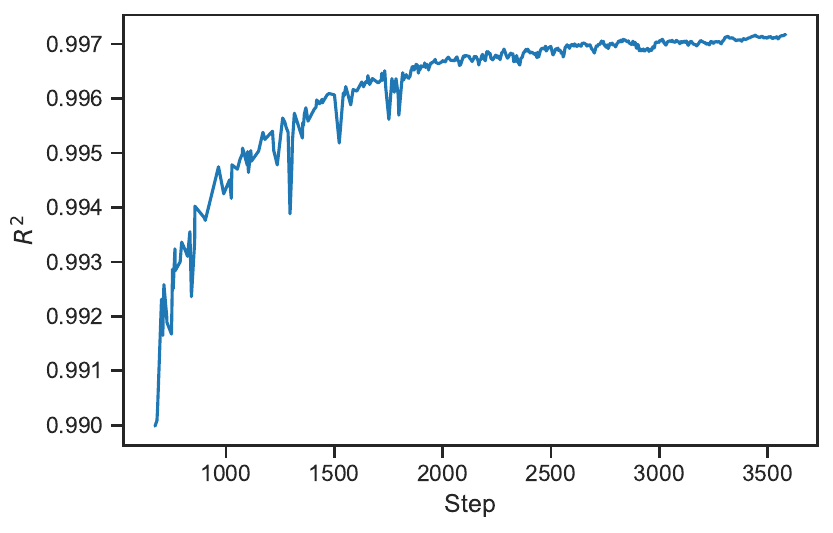}

}

\subcaption{\label{fig-val-2}\(R^{2}\)}

\end{minipage}%

\caption{\label{fig-val}Evaluation metrics}

\end{figure}

The metrics for validation were the mean absolute error (MAE) and
\(R^2\). Last values are shown in Table~\ref{tbl-eval}, , while the
evolution along the various steps is in Figure~\ref{fig-val}.

\begin{figure}

\centering{

\includegraphics{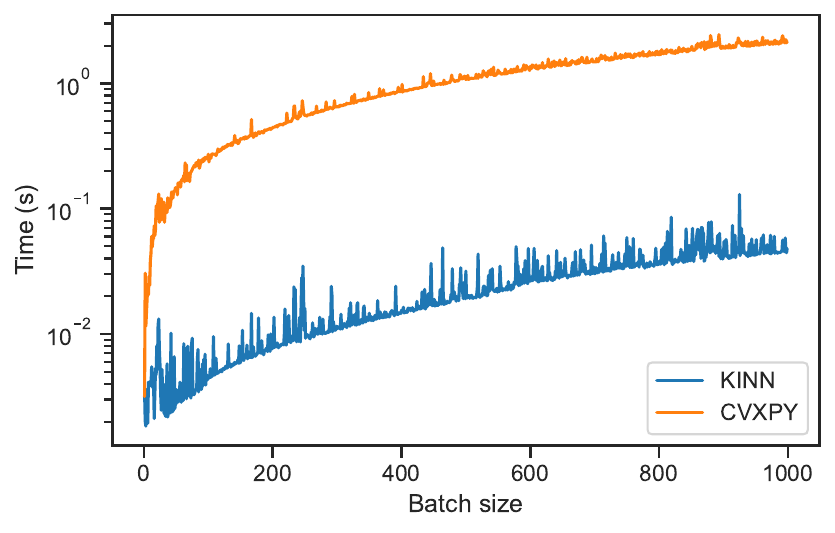}

}

\caption{\label{fig-times}Computation time comparison}

\end{figure}

By increasing the number of points, an inference time comparison was
performed on an Apple M2 Pro processor with backend for PyTorch's MPS.
The difference is remarkable, about two orders of magnitude
(Figure~\ref{fig-times}): with a batch size of \(1000\), cvxpylayers
took \(2.35~\textrm{s}\) as opposed to KINN's \(0.06~\textrm{s}\).

\begin{longtable}[]{@{}cc@{}}
\caption{Evaluation metrics}\label{tbl-eval}\tabularnewline
\toprule\noalign{}
Mertric & Value \\
\midrule\noalign{}
\endfirsthead
\toprule\noalign{}
Mertric & Value \\
\midrule\noalign{}
\endhead
\bottomrule\noalign{}
\endlastfoot
MAE & 0.0056 p.u. \\
\(R^{2}\) & 0.9972 \\
\end{longtable}

\section{Conclusions}\label{conclusions}

KKT-Informed Neural Network (KINN) was introduced as a neural
network-based approach for solving parametric convex optimization
problems. The method leverages the Karush-Kuhn-Tucker (KKT) conditions
as a learning bias, integrating them into the loss function to ensure
that the network's predictions adhere to the necessary optimality
criteria. This allows the network to efficiently estimate solutions
while sacrificing some degree of guaranteed optimality in favor of
significant improvements in computational speed and scalability.

The experimental results from the provided test case demonstrated that
KINN is highly effective in providing near-optimal solutions while
enabling parallel problem-solving. The comparison with traditional
tools, highlighted KINN's ability to solve optimization problems in a
fraction of the time, with minimal loss of accuracy. Metrics such as
mean absolute error (MAE) and \(R^2\) confirm that KINN produces
reliable solutions within acceptable tolerances for real-world
applications. Although some optimality is traded for speed, KINN
provides a highly practical tool for scenarios where rapid
decision-making is critical.

In future work, the potential for expanding this approach to handle
non-convex optimization problems will be explored. Additionally, hybrid
architectures that combine KINN with traditional optimization methods
may offer further improvements in performance, providing a balance
between speed and guaranteed optimality across a broader range of
problem domains.

  \bibliography{Preprint.bib}

\end{document}